\documentclass[review,11pt]{elsarticle}

\usepackage{amssymb}
\usepackage{xcolor}
\usepackage{amsmath}
\usepackage{amsthm}
\usepackage[flushleft]{threeparttable} 
\usepackage{subcaption}
\usepackage{algorithm}
\usepackage{algpseudocode}

\journal{Mathematical Biosciences}

\newcommand{\fder}[2]{\dfrac{\mathrm{d}#1}{\mathrm{d}#2}}

\begin{document}

\begin{frontmatter}



\title{Persistent instability in a nonhomogeneous delay differential equation system of the Valsalva maneuver}


\author{E. Benjamin Randall, Nicholas Z. Randolph, and Mette S. Olufsen}

\address{Department of Mathematics, North Carolina State University, Raleigh, NC}

\begin{abstract}
Delay differential equations are widely used in mathematical modeling to describe physical and biological systems, often inducing oscillatory behavior. In physiological systems, this instability may signify (i) an attempt to return to homeostasis or (ii) system dysfunction. In this study, we analyze a nonlinear, nonautonomous, nonhomogeneous open-loop neurological control model describing the autonomic nervous system response to the Valsalva maneuver (VM). We reduce this model from 5 to 2 states (predicting sympathetic tone and heart rate) and categorize the stability properties of the reduced model using a two-parameter bifurcation analysis of the sympathetic delay ($D_s$) and time-scale ($\tau_s$). Stability regions in the $D_s \ \tau_s$-plane for this nonhomogeneous system and its homogeneous analog are classified numerically and analytically identifying transcritical and Hopf bifurcations. Results show that the Hopf bifurcation remains for both the homogeneous and nonhomogeneous systems, while the nonhomogeneous system stabilizes the transition at the transcritical bifurcation. This analysis was compared with results from blood pressure and heart rate data from three subjects performing the VM: a control subject exhibiting sink behavior, a control subject exhibiting stable focus behavior, and a patient with postural orthostatic tachycardia syndrome (POTS) also exhibiting stable focus behavior. Results suggest that instability caused from overactive sympathetic signaling may result in autonomic dysfunction.
\end{abstract}




\end{frontmatter}



\section{Introduction}
\label{Intro}

\par \noindent Analysis of the dynamics of differential equations models can shed light on model prediction outcomes reflecting healthy versus disease states. Model predictions of normal and abnormal behavior often require either changes in the nominal parameter values or the modeled dynamic pathways. For the former, healthy model predictions can be a result of operating in a stable region of the parameter space, while diseased outcomes can be a result of (i) a change in parameter regime to an unstable region or (ii) a bifurcation to a new stable equilibrium. To explore this further, we conduct a stability analysis to determine where in the parameter space instability arises in a mathematical model describing the autonomic nervous response to the Valsalva maneuver (VM). 

\par The model used in this study is formulated using a system of nonlinear delay differential equations (DDEs), which are common in the study of many real-world systems \cite{Banks2017,Keane2017,Nelson2002,Vielle1998,Wilhelm2009}.  DDEs are known to change the dynamical behavior, causing bistability or instability in some systems \cite{dOnofrio2010,Bi2013,Rao2018,Sipahi2011} and broadening the stable region of others \cite{Sipahi2011}. Physically, delays are often used to avoid adding equations describing the process causing the delay. In our case, the process is the transmission of the sympathetic response along the sympathetic ganglia chain. Given DDEs are known to generate instability \cite{Sipahi2011}, it is important to analyze whether a delay is critical to model the system and to test alternative formulations, such as distributed delays, which impose chains of differential equations with varying time-scales \cite{Olsen2014}.  Often, distributed delays garner a similar effect as a discrete DDE without the added computational expenditure but at the cost of increasing the dimension of the state and parameter spaces. The choice to use distributed versus discrete delays is problem-dependent. To avoid an increase in the dimension of the state space, we model the effect of the sympathetic control using a discrete time-delay. 

\par Numerical tools for bifurcation analysis of DDEs exist, {\it e.g.}, \texttt{DDE-Biftool} \cite{Engelborghs2002} and \texttt{knut} \cite{Roose2007,Szalai2006}. The former is a powerful collection of MATLAB$^\text{\textregistered}$ routines for autonomous DDEs with constant and state-dependent delays, which has previously been used in two-parameter bifurcation analyses \cite{Krauskopf2014,Luzyanina2005,Makroglou2006}. However, in this study, we analyze a nonautonomous system of stiff DDEs, which \texttt{DDE-Biftool} currently cannot accommodate. \texttt{DDE-Biftool} uses the built-in MATLAB delay solver \texttt{dde23}, which does not account for stiff systems with multiple time-scales. \texttt{knut} \cite{Szalai2006,Roose2007} is a bifurcation analysis package in C$^{\text{++}}$ that allows periodic forcing functions. However, the forcing term for this model is not periodic. In this study, we evaluate the forward model over a discretized mesh of the parameter space and categorize the model output into groups analogous to the harmonic oscillator: critically damped sink, overdamped sink, stable focus, limit cycle, and unstable. 

\par Unstable modes arise in many physical and biological systems naturally and avoiding these modes is of particular interest in recent years \cite{Bi2013,Cheng2018,Rao2018,Wilhelm2009}. For physiological processes, at rest the body is mainly operating via negative feedback mechanisms that maintain homeostasis, {\it e.g.}, the baroreceptor reflex ({\it baroreflex}) modulating blood pressure and heart rate. However, it is known that in some disease states the negative feedback mechanisms fail and are overridden by positive feedback mechanisms, {\it e.g.}, the Bezold-Jarisch reflex during vasovagal syncope, which causes the system to transition to an unstable state (syncope) \cite{Hall2016}. In this study, the objective is to characterize the stability regions important to autonomic dysfunction (AD) in patient data. More specifically, we investigate the persistent instability as a result of the baroreflex response to the VM. 

\par
\begin{figure}[t!]
	\centering 
	\includegraphics[width=\textwidth]{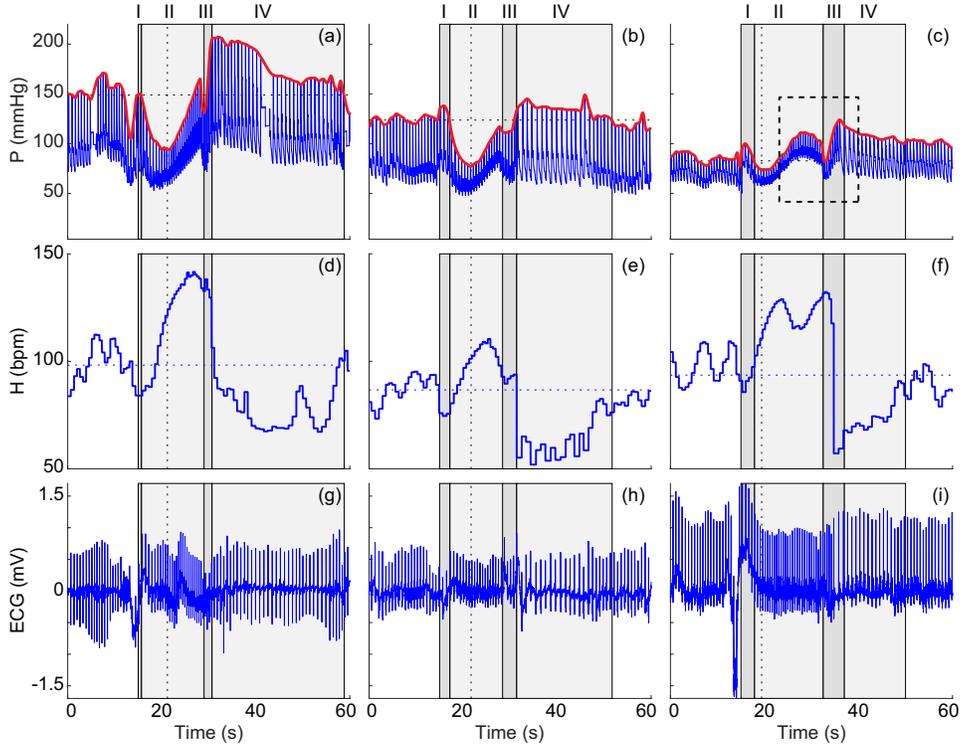}
	\caption{Blood pressure ($P$, mmHg) with systolic blood pressure (red) indicated, heart rate ($H$, bpm), and electrocardiogram (ECG, mV) data for each subject. Valsalva maneuver phases are indicated with alternating gray (I and III) and light gray (II and IV) boxes. Early and late phase II is divided with a vertical dashed line. (a, d, g) Subject 1 - control subject with sink behavior. (b, e, h) Subject 2 - control subject with stable focus behavior. (c, f, i) Subject 3 - patient with postural orthostatic tachycardia syndrome (POTS) exhibiting M behavior (dashed black box) with stable focus behavior. Descriptions for subjects are given in Table \ref{datatable}.} 
	\label{data}
\end{figure}

 We aim to categorize disease and healthy states based on a two-parameter bifurcation analysis. We show that the sympathetic time-scale, $\tau_s$, and sympathetic delay, $D_s$, parameters are intrinsically linked and interactions between them can cause oscillations and unstable behavior. The use of stability analysis to examine effects of the delay in the baroreflex response has been done in one study by Ottesen \cite{Ottesen1997}. He  performed a two-parameter bifurcation analysis and showed that when the time delay is varied over its physiological range, stability switches arise. However, oscillatory modes were not compared to patient data. In this study, we compare the sympathetic outflow and heart rate responses of two control subjects and a patient with AD exhibiting the M response to the VM as categorized by Palamarchuk {\it et al.} \cite{Palamarchuk2016} (Figure \ref{data}), determining parameter regimes where instability occurs. The M behavior is hypothesized as overactive sympathetic and parasympathetic activity. In this study, we analytically determine the locations of the different stability regimes by solving and comparing the homogeneous DDE system to numerical simulations of the associated nonhomogeneous system and discuss the effects of the forcing function on the stability. In addition, we connect these results to physiological data exhibiting instability.

\section{Materials and methods}
\label{Modeldevelopment}

\par \noindent This study analyzes a neurological control model of the autonomic nervous response to the VM. Utilizing systolic blood pressure (SBP, mmHg) and thoracic pressure ($P_{th}$, mmHg) as inputs, the model predicts heart rate and sympathetic and parasympathetic nervous system responses. 

\begin{table}[!t]
	\centering 
	\footnotesize
	\begin{threeparttable} 
		\caption{Subject data.}
		\begin{tabular}{clclcc}
			\hline
			{\bf Subject} & \multicolumn{1}{c}{{\bf Classification}}  & {\bf Age} & \multicolumn{1}{c}{{\bf Sex}} & $\bar{P}$   &  $\bar{H}$ \\
			& &  (years) & & (mmHg)  & (bpm)  \\ 
			\hline 
			\hline 
			1 & Control & 21 & Female &  149  & 98 \\ 
			2 & Control & 27 & Male  &  117  &  87 \\
			3 & POTS with M behavior  & 57 & Female &  83  &  94  \\
			\hline 
		\end{tabular} 
		\label{datatable}
		\begin{tablenotes} 
			\small
			\item POTS - postural orthostatic tachycardia syndrome. 
			\item  $\bar{P}$ - mean systolic blood pressure. 
			\item  $\bar{H}$ - mean heart rate. 
		\end{tablenotes} 
	\end{threeparttable} 
\end{table} 

\subsection{Data}

\par \noindent Blood pressure ($P$, mmHg) and electrocardiogram (ECG, mV) measurements were collected via a Finometer (Finapres Medical Systems BV, Amsterdam, The Netherlands) and a precordial ECG-lead, respectively, and saved in LabChart$^\text{\textregistered}$ for three subjects performing a Valsalva maneuver. All subjects gave consent to participate in this study  and the protocol was approved by the Ethics Committee for the Capital Region, Denmark. Heart rate ($H$, bpm) was computed in LabChart$^\text{\textregistered}$ using cyclic detection for human ECG. Table \ref{datatable} summarizes the data used in this study. The blood pressure, heart rate, and ECG data for all subjects are shown in Figure \ref{data}. SBP is calculated as the interpolation of consecutive local maxima in the blood pressure (Figures \ref{data}a-c red curve). Subjects 1 and 2 exhibit no AD  and Subject 3 has postural orthostatic tachycardia syndrome (POTS), determined as an increase in heart rate of $\geq$30 bpm without an associated increase in blood pressure during a postural change \cite{Weimer2010}. Subject 3 exhibits the M blood pressure response to the VM as categorized by Palamarchuk {\it et al.} \cite{Palamarchuk2016} (Figure \ref{data}c dashed  black box). 

\subsection{Valsalva maneuver}

\par \noindent The Valsalva maneuver (VM) is a clinical test that involves forced expiration while maintaining an open glottis \cite{Hamilton1943}. In response to a sudden decrease in blood pressure, the VM initiates the baroreflex, which inhibits parasympathetic and stimulates sympathetic activity, increasing heart rate \cite{Boron2017}.The VM is divided into four phases (illustrated in Figure \ref{data}):
\begin{enumerate}[I.]
	\item The breath hold causes a sharp increase in blood pressure and slight decrease in heart rate.
	\item Phase II is divided into two sections:
	\begin{enumerate}[i.]
		\item Early phase II: Blood pressure drops below baseline significantly, triggering parasympathetic withdrawal and heart rate acceleration.
		\item Late phase II: Delayed sympathetic activation accelerates heart rate further and increases peripheral vascular resistance, resulting in an increase in blood pressure. 
	\end{enumerate}  
	\item Release of the breath hold causes a sharp decrease in blood pressure, triggering a second parasympathetic withdrawal. 
	\item Increased sympathetic activation causes blood pressure to overshoot and return to baseline within 30 s, while normalization of parasympathetic activity causes a sharp drop in heart rate and subsequent return to baseline. 
\end{enumerate} 

\begin{figure}[t]
	\centering 
	\includegraphics[width=\textwidth]{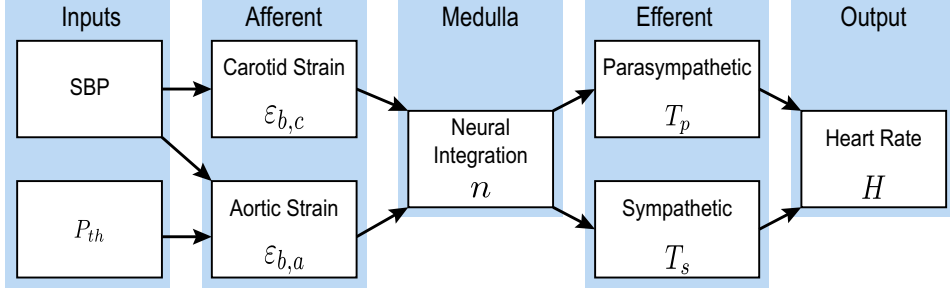}
	\caption{Baroreflex model schematic. Systolic blood pressure (SBP, mmHg) and thoracic pressure ($P_{th}$, mmHg) are inputs to the model. The afferent signals are propagated via the carotid baroreceptor strain ($\varepsilon_{b,c}$, dimensionless) stimulated by SBP and the aortic baroreceptor strain ($\varepsilon_{b,a}$, dimensionless) stimulated by the difference in SBP and $P_{th}$. These signals are integrated in the medulla ($n$, s$^{-1}$), initiating an efferent signal cascade via the parasympathetic ($T_p$, dimensionless) and sympathetic ($T_s$, dimensionless) nervous systems, which modulates heart rate ($H$, bpm).} 
	\label{modelschematic} 
\end{figure}

\begin{table}[!ht]
	\centering 
	\footnotesize  
	\renewcommand{\arraystretch}{1.5}
	\begin{threeparttable} 
		\caption{Initial conditions and constant history value.}
		
		\begin{tabular}{lrl}
			\hline 
			\multicolumn{1}{c}{{\bf Description}} & \multicolumn{2}{c}{{\bf Value}} \\ 
			\hline
			Carotid baroreceptor strain & $\varepsilon_{b,c,0} =$ & $  1 - \sqrt{\dfrac{2}{A + 1}}$ \\
			Aortic baroreceptor strain & $\varepsilon_{a,c,0} =$ & $  1 - \sqrt{\dfrac{2}{A + 1}}$ \\
			Parasympathetic outflow & $T_{p,0} =$ & 0.8 \\ 
			Sympathetic outflow & $T_{s,0} =$ & 0.2 \\ 
			Heart rate & $H_{0} = $ & $\bar{H}$ \\
			\hline 
		\end{tabular} 
		\label{init}
		\begin{tablenotes} 
			\small
			\item $\bar{H}$ is the baseline heart rate given in Table \ref{datatable}. 
			\item $A$ is given in Table \ref{parameters}. 
		\end{tablenotes} 
	\end{threeparttable} 
	\renewcommand{\arraystretch}{1} 
\end{table} 

\subsection{Model overview} 

\par \noindent The model analyzed in this study is from our previous work in Randall {\it et al}. \cite{Randall2019} but does not account for the respiratory sinus arrhythmia (RSA). Since the RSA sub-model does not depend on the effect the sympathetic delay, we remove this model component and focus solely on the baroreflex sub-model. A schematic of the baroreflex model is given in Figure \ref{modelschematic}. 

\par
\begin{figure}[!t]
	\centering 
	\includegraphics[]{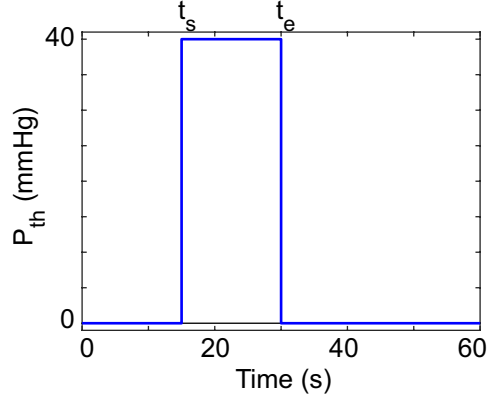}
	\caption{Thoracic pressure ($P_{th}$, mmHg) to induce the Valsalva maneuver calculated in equation \eqref{pth}.} 
	\label{Pthfig} 
\end{figure}

 The open-loop model takes SBP data and $P_{th}$ as inputs predicting heart rate and parasympathetic and sympathetic responses. For simplicity, the thoracic pressure is modeled as 
\begin{equation} 
P_{th}(t) = \left\{
\begin{array}{ll} 
40 & \text{for} \ t_s \leq t \leq t_e \\ 
0 & \text{otherwise},
\end{array}
\right. 
\label{pth} 
\end{equation} 
where $t_s$ and $t_e$ are the start and end of the VM from the data as shown in Figure \ref{Pthfig}. Due to the delay in sympathetic signal transduction, the model incorporates a discrete time-delay, $D_s$, into the  differential equation modeling sympathetic outflow, giving the following system of equations: 
\begin{equation}
\renewcommand{\arraystretch}{1.5}
\left\{
\begin{array}{rll}
\fder{\varepsilon_{b,c}}{t} &= \dfrac{-\varepsilon_{b,c}+K_b\varepsilon_{w,c}}{\tau_b}, & \varepsilon_{b,c}(0) = \varepsilon_{b,c,0},\\ 
\fder{\varepsilon_{b,a}}{t} &= \dfrac{-\varepsilon_{b,a}+K_b\varepsilon_{w,a}}{\tau_b}, & \varepsilon_{b,a}(0) = \varepsilon_{b,a,0}, \\
\fder{T_{p}}{t} &= \dfrac{-T_{p} + K_{p}G_{p}}{\tau_{p}}, & T_{p}(0) = T_{p,0}, \\ 
\fder{T_{s}}{t} &= \dfrac{-T_{s}(t-D_s) + K_{s}G_{s}}{\tau_{s}}, & T_{s}(t) = T_{s,0}, \ t\in[-D_s,0], \\ 
\fder{H}{t} &= \dfrac{-H+\tilde{H}}{\tau_H}, & H(0) = H_0,
\end{array}
\right.
\renewcommand{\arraystretch}{1}
\label{fullsystem}
\end{equation}
where $\tilde {H} = H_I(1-H_{p}T_{p}+H_s T_s)$. The states are the carotid ($\varepsilon_{b,c}$, dimensionless) and aortic ($\varepsilon_{b,a}$, dimensionless) baroreceptor strains, the parasympathetic ($T_p$, dimensionless) and sympathetic ($T_s$, dimensionless) outflows, and the heart rate ($H$, bpm).  $K_l$ (dimensionless) and $\tau_l$ (s)  for $l = b$, $p$, $s$, or $H$ denoting baroreceptor, parasympathetic, sympathetic, and heart rate, respectively, are the gains and time-scales for each of the differential equations with units and values given in Table \ref{parameters}. $H_I$ (bpm) denotes the intrinsic heart rate of 100 bpm, and $H_p$ (dimensionless) and $H_s$ (dimensionless) are gains scaling parasympathetic and sympathetic outflow, respectively. Initial conditions and constant history value are summarized in Table \ref{init}. Arterial wall strain $\varepsilon_{w,j}$ for $j = c$ or $a$ denoting the carotid sinus and aortic arch, respectively, is a nonlinear sigmoid-like function predicting arterial wall deformation given by 
\begin{equation}
\varepsilon_{w,j} = 1 - \sqrt{\dfrac{1+e^{-q_w(P_j-s_w)}}{A+e^{-q_w(P_j-s_w)}}}
\label{epswj}
\end{equation}
for the carotid ($P_c = \text{SBP}$) and aortic ($P_a = \text{SBP} - P_{th}$) pressures, where $q_w$ (mmHg$^{-1}$) and $s_w$ (mmHg) are the steepness and half-saturation values and $A$ (dimensionless) is an offset parameter. The saturation functions $G_l$ for $l = p$ or $s$ are the sigmoidal relations
\begin{equation}
\begin{array}{rcl}
G_{p} &=& \dfrac{1}{1+e^{-q_{p}(n-s_{p})}} \quad \text{and} \\ 
G_s &=& \dfrac{1}{1+e^{q_s(n-s_s)}},
\end{array}
\end{equation}
where $q_l$ (s) and $s_l$ (s$^{-1}$) are the steepness and half-saturation values and 
\begin{equation} 
n = B(\varepsilon_{w,c}-\varepsilon_{b,c}) + (1-B)(\varepsilon_{w,a}-\varepsilon_{b,a}), \quad B \in[0,1]
\label{n}
\end{equation}  
is a convex combination of the relative strains.

\par In summary, the model consists of 5 differential state equations and 20 parameters with one DDE. The model has the form 
\begin{equation}
\fder{\mathbf{x}}{t}(t) = f(t, \mathbf{x}(t), \mathbf{x}(t-D_s); \mathbf{\theta}), \quad \mathbf{x}(t) = \mathbf{x}_0, \ t\in [-D_s, 0],	
\label{ddeform}
\end{equation} 
where $f$ is the right hand side, $\mathbf{x} = [\varepsilon_{b,c}, \varepsilon_{b,a}, T_{p}, T_s, H]^T \in \mathbb{R}^5$ is the state vector, $\mathbf{x}_0 = [\varepsilon_{b,c,0}, \varepsilon_{b,a,0}, T_{p,0}, T_{s,0}, \bar{H}]^T\in \mathbb{R}^5$ is the constant history vector, $D_s$ is the discrete delay, and $\theta \in \mathbb{R}^{20}$ is the parameter vector
\begin{equation}
\theta = [ A, B, K_b, K_{p}, K_s, \tau_b, \tau_{p}, \tau_s, \tau_H, q_w, q_{p}, q_s, s_w, s_{p}, s_s, H_I, H_{p}, H_s, t_s, t_e]^T,
\end{equation}
where each $\theta_i > 0$. Nominal parameter values are summarized in Table \ref{parameters} and an explanation of parameter assignments is given in the Appendix. 

\begin{table}[p]
	\centering 
	\footnotesize  
	\begin{threeparttable}
		\caption{Nominal parameter values.}
		\begin{tabular}{cclc} 
			\hline 
			{\bf Symbol} & {\bf Units} & \multicolumn{1}{c}{{\bf Description}} & {\bf Nominal}\\
			& & & {\bf value} \cite{Randall2019}\\
			\hline 
			\hline
			$A$ &  & Cross-sectional area ratio & 5  \\
			$B$ & s$^{\text{-1}}$ & Neural integration parameter & 0.5 \\ 
			$K_b$ &  & Baroreceptor strain gain & 0.1 \\ 
			$K_{p}$ & & Baroreflex-parasympathetic gain & 5 \\
			$K_{s}$ & & Baroreflex-sympathetic gain & 5 \\ 
			$\tau_b$ & s & Baroreceptor strain time-scale & 0.9  \\ 
			$\tau_{p}$ & s & Baroreflex-parasympathetic time-scale & 1.8  \\ 
			$\tau_s$ & s & Baroreflex-sympathetic time-scale & 10  \\ 
			$\tau_H$ & s & Heart rate time-scale & 0.5 \\
			$q_w$ & mmHg$^{\text{-1}}$ & Arterial wall strain sigmoid steepness & 0.04  \\ 
			$q_{p}$ & s & Baroreflex-parasympathetic sigmoid steepness & 10  \\ 
			$q_s$ & s & Baroreflex-sympathetic sigmoid steepness & 10 \\ 
			$s_w$ & mmHg & Arterial wall strain half-saturation value  & $\bar{P}$\\ 
			$s_{p}$ & s$^{\text{-1}}$ & Baroreflex-parasympathetic half-saturation value & 0.55 \eqref{s_p}  \\ 
			$s_{s}$ & s$^{\text{-1}}$ & Baroreflex-sympathetic half-saturation value & 0.06 \eqref{s_s}   \\  
			$H_I$ & bpm  & Intrinsic heart rate&  100  \\ 
			$H_{p}$ & & Baroreflex-parasympathetic heart rate gain & 0.22$\pm$0.04 \eqref{H_p}    \\ 
			$H_{s}$ & & Baroreflex-sympathetic heart rate gain&  0.37$\pm$0.2 \eqref{H_s}   \\
			$D_s$ & s & Sympathetic delay parameter & 3  \\
			$t_s$ & s & Valsalva maneuver start time & data \\ 
			$t_e$ & s & Valsalva maneuver end time & data \\
			\hline
		\end{tabular}  
		\label{parameters} 
		\begin{tablenotes}
			\small
			\item A blank space in the Units column indicates the parameter is dimensionless.
			\item A mean $\pm$ one SD is reported where applicable. 
			\item $\bar{P}$ denotes the baseline systolic blood pressure as listed in Table \ref{datatable}. 
			\item ``data" refers to the fact that this quantity was taken directly from the data. 
		\end{tablenotes} 	
	\end{threeparttable}
\end{table}  

\subsection{Physiological model reduction}

\par \noindent To analyze of the effect and dynamics of the DDE, we reduce the model to two states, the baroreflex-mediated sympathetic tone $T_s$ (the DDE) and heart rate $H$, which is affected by the DDE. Since the model is open-loop, the other three states ($\varepsilon_{b,c}$, $\varepsilon_{b,a}$, and $T_p$) are not impacted by the DDE. To establish the two-state model, we eliminate these states by forming algebraic relations, taking advantage of short time-scales as explained below. 

\par First, since the model is open-loop, there is no feedback to the states that come before the delayed state $T_s$, {\it i.e.}, $\varepsilon_{b,c}$, $\varepsilon_{b,a}$, and $T_p$,  as shown in Figure \ref{modelschematic}. For the purposes of this analysis, we assume these states remain in steady-state for the entire time interval, and hence, we can reduce them algebraically. We reformulate the differential equation for $T_p$ as 
\begin{equation}
\begin{array}{rcl} 
\fder{T_{p}}{t} &=& \dfrac{-T_{p} + K_{p} G_{p}}{\tau_{p}} \\
 \Rightarrow \quad \tau_{p} \fder{T_{p}}{t} &=& -T_{p} + K_{p} G_{p}.
\end{array}
\end{equation} 
Since the time-scale $\tau_{p}$ is an order of magnitude smaller than $\tau_s$, we remove this differential equation by setting $\tau_{p} = 0$ and solving for $T_{p}$, giving
\begin{equation} 
T_{p}= K_{p} G_{p}. 
\end{equation} 
We make a similar simplification for the baroreceptor strains ($\varepsilon_{b,j}$ for  $j = c$ or $a$), taking advantage of the small time-scale $\tau_b$ in relation to the magnitude of $\tau_s$. Thus, 
\begin{equation} 
\varepsilon_{b,j} = K_b \varepsilon_{w,j}.
\end{equation} 

\par
\begin{figure}[!t]
	\centering 
	\includegraphics[]{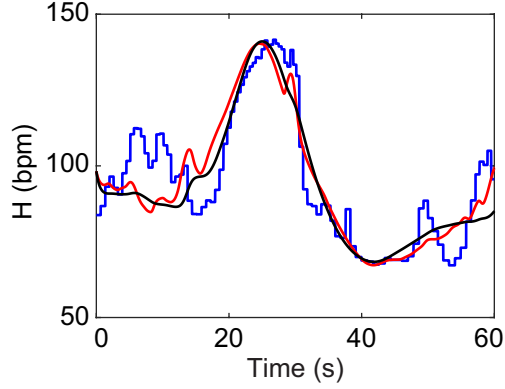}
	\caption{Full five-state (red) and reduced two-state (black) model fits to heart rate data (blue) for Subject 1.}
	\label{fullvstwostate}
\end{figure}

 Second, since we are interested in the instability caused by the VM, we accentuate the effect of the increased thoracic pressure ($P_{th}$) during to the VM by eliminating the carotid pathway. With this simplification, the model depends solely on the effect of the aortic baroreceptors, that is, we set $B = 0$ and 
\begin{equation} 
n = \varepsilon_{w,a} - \varepsilon_{b,a} = (1 - K_b) \varepsilon_{w,a} = (1 - K_b) \Bigg(1 - \sqrt{ \dfrac{ 1 + e^{-q_w (P_a - s_w)}}{A + e^{-q_w (P_a - s_w)}}} \Bigg). 
\end{equation} 
The resulting model is a system of two states, $T_s(t)$ and $H(t)$, of the form  
\begin{equation}
\renewcommand{\arraystretch}{1.5}
\left\{
\begin{array}{rll}
\fder{T_s}{t} &= \dfrac{-T_s(t-D_s) + K_sG_s}{\tau_s}, & T_s(t) = T_{s,0}, t\in[-D_s,0],\\
\fder{H}{t} &= \dfrac{-H+H_I (1 - H_{p} K_{p} G_{p} + H_s T_s)}{\tau_H}, & H(0) = H_0.
\end{array}
\right.
\renewcommand{\arraystretch}{1.5}
\label{reducedsystem}
\end{equation}
The reduced system is in the form of equation \eqref{ddeform}, where $\mathbf{x} = [T_s, H]^T$ $\in \mathbb{R}^2$, $\mathbf{x}_0 = [T_{s,0}, H_0]^T$ $\in \mathbb{R}^2$ for $t \in [-D_s, 0]$, $D_s$ is a discrete delay, and $\theta \in \mathbb{R}^{16}$ is the vector of parameters 
\begin{equation}
\theta = [A, K_{p}, K_s, \tau_s, \tau_H, q_w, q_{p}, q_s, s_w, s_{p}, s_s, H_I, H_{p}, H_s, t_s, t_e]^T.
\end{equation}
The model outputs for the full model (red, equation \eqref{fullsystem}) and the reduced model (black, equation \eqref{reducedsystem}) in comparison to the heart rate data from Subject 1 are shown in Figure \ref{fullvstwostate}. 

\par
\begin{figure}[!t]
	\centering 
	\includegraphics[width=\textwidth]{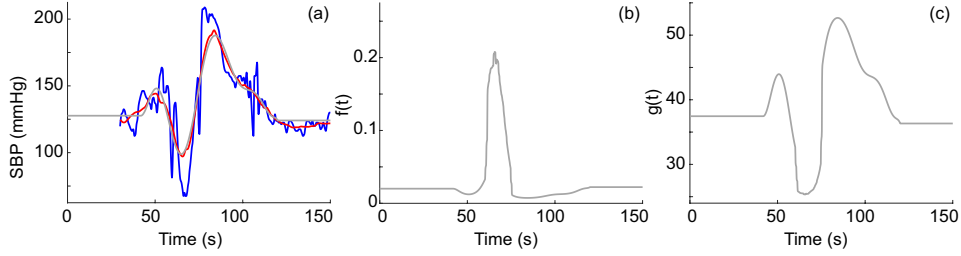}
	\caption{Forcing functions. (a) Systolic blood pressure (SBP, mmHg) data from Subject 1 (blue), filtered data using the \texttt{movmean} command in MATLAB$^\text{\textregistered}$ with a window of one second (red), and a 10$^{\text{th}}$ degree polynomial fit to the \texttt{movmean} signal (gray). (b) $f(t)$ (s$^{-1}$) for $T_s$ given in equation \eqref{f}. (c) $g(t)$ (bpm s$^{-1}$) for $H$ given in equation \eqref{g}.}
	\label{forcingfcns}
\end{figure}

 The model in equation \eqref{reducedsystem} is linear with respect to the states and reformulating the system yields
\begin{align}
\fder{}{t} \left[
\begin{array}{c}
T_s \\ 
H 
\end{array} 
\right] &= 
\left[
\begin{array}{l}
-\dfrac{1}{\tau_s}T_s(t-D_s) + f(t) \vspace*{0.75mm} \\ 
-\dfrac{1}{\tau_H}H+ \dfrac{H_IH_s}{\tau_H}T_s + g(t)
\end{array}
\right], \\ 
\left[
\begin{array}{c}
T_s(t) \\ 
H(t) 
\end{array} 
\right] &= 
\left[
\begin{array}{l}
T_{s,0} \\ 
H_0
\end{array}
\right],
 t\in[-D_s, 0].
\label{system1}
\end{align}
The forcing functions $f$ (s$^{-1}$) and $g$ (bpm s$^{-1}$) for $T_s$ and $H$, respectively, are given by
\begin{equation}
f(t) = \dfrac{K_s}{\tau_s} G_s \Big( n \big( P_a(t) \big) \Big) 
\label{f}
\end{equation}
and 
\begin{equation}
g(t) = \dfrac{H_I}{\tau_H} \bigg( 1 - H_{p} K_{p} G_{p} \Big( n \big( P_a(t) \big) \Big)  \bigg).
\label{g}
\end{equation}
This is a nonhomogeneous, nonautonomous DDE system. Forcing functions $f(t)$ and $g(t)$ represent the dynamics induced by the blood pressure responses to the VM. Since the forcing functions use blood pressure data as an input, we ensure smoothness by filtering the data using the \texttt{movmean} command in MATLAB$^\text{\textregistered}$ with a window of one second. Then, we fit a $10^\text{th}$ degree polynomial to the filtered signal, that is, the coefficients $a_i$ of a polynomial $P = \sum\limits_{i=0}^{10} a_ix^i$ were optimized to fit the filtered signal. Polynomials of orders $>$10 produced high frequency oscillatory behavior at baseline. We artificially extended the SBP before the curve to ensure the model began in steady-state and after to accentuate the oscillatory behavior of the signal if it arises. Figure \ref{forcingfcns_a} shows the original blood pressure data (blue), the moving mean (red), and the fitted polynomial $P$ (gray). Figures \ref{forcingfcns}b and c display the forcing functions $f$ and $g$, respectively. 

\par In summary, we consider the system
\begin{equation}
\fder{\mathbf{x}}{t} = \mathbf{A}\mathbf{x} + \mathbf{B}\mathbf{x}_{D_s} + \mathbf{f}, \quad \mathbf{x}(t) = \mathbf{x}_0, t\in[-D_s,0],
\label{matsys}
\end{equation}
where $\mathbf{x}_{D_s} = \mathbf{x}(t - D_s)$ is the vector of delayed states, $\mathbf{A}$ and $\mathbf{B}$ are constant matrices given as
\begin{equation}
\mathbf{A} = \begin{bmatrix}
0 & 0 \\ 
\dfrac{H_I H_s}{\tau_H} & -\dfrac{1}{\tau_H}
\end{bmatrix}
\quad 
\text{and} 
\quad 
\mathbf{B} = \begin{bmatrix}
-\dfrac{1}{\tau_s} & 0 \\ 
0 & 0
\end{bmatrix},
\end{equation}
and $\mathbf{f}$ is the forcing vector given by
\begin{equation}
\mathbf{f} = \left[ 
\begin{array}{c}
f(t) \\ 
g(t) 
\end{array}
\right].
\label{forcingvector} 
\end{equation} 

\section{Stability analysis}

\par \noindent The stability of the DDE in equation \eqref{matsys} depends on both the homogeneous solution and the effect of the forcing function. In this section, we analytically explore the homogeneous equation by classifying the roots of the CE, $\phi$. We numerically categorize the behavior of the nonhomogeneous, nonautonomous system using Algorithm \ref{alg1} that takes advantage of the gradient of the solution after the VM occurs.  

\subsection{Homogeneous system}

\par \noindent To analyze the stability of equation \eqref{matsys}, we first consider the homogeneous equation 
\begin{equation}
\fder{\mathbf{x}}{t} = \mathbf{A}\mathbf{x} + \mathbf{B}\mathbf{x}_{D_s}, \quad \mathbf{x}(t) = \mathbf{x}_0,  t\in[-D_s,0], 
\label{homogeneous}
\end{equation}
which is expanded as 
\begin{align}
\fder{T_s}{t} &= -\frac{1}{\tau_s}T_s(t-D_s), & & T_s(t) = T_{s,0}, t\in[-D_s,0] \label{homodT}, \\
\fder{H}{t} &= \frac{1}{\tau_H}\big(-H+ H_IH_sT_s\big), & &  H(0) = H_0. 
\label{homodH}
\end{align}
In this system, the origin is a unique critical point, {\it i.e.}, if another critical point exists, then from equation \eqref{homodH}, we have
\begin{equation} 
0 = \dfrac{1}{\tau_H} (-H + H_I H_s T_s) \quad \Rightarrow \quad H = H_I H_s T_s. 
\end{equation}  
From \eqref{homodT}, we have 
\begin{equation} 
0 = \dfrac{1}{\tau_s} T_s(t - D_s) \quad \Rightarrow \quad T_s(t) = 0. 
\end{equation}
Therefore, the origin is the only critical point of the system in \eqref{homogeneous}. 

\par
\begin{figure}[!t]
	\centering 
	\includegraphics[width=.8\textwidth]{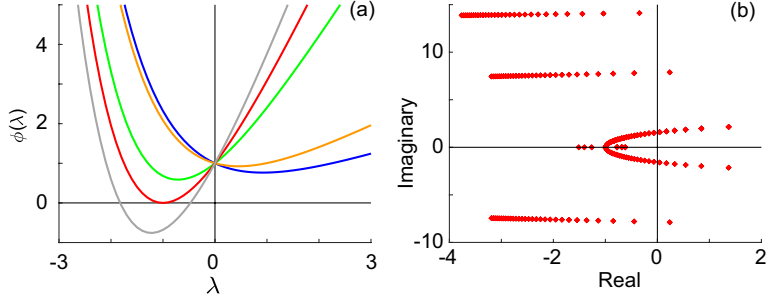}
	\caption{(a) The characteristic equation (CE), solving $\phi(\lambda) = 0$ (given in equation \eqref{chareq}) for $D_s = 1$ and $\tau_s$ varied. The graph displays the CE of each type: gray curve ($\tau_s =  1.25e$), the CE has two unique real solutions ($\lambda_1\neq \lambda_2 \in \mathbb{R}$); red curve ($\tau_s = e$), the CE has one real solution ($\lambda_1 = \lambda_2 = \lambda \in \mathbb{R}$); green, orange, and blue curves ($\tau_s = 0.75e$, $2/\pi$, and $1.25/\pi$, respectively), the CE has infinitely many solutions ($\lambda = \alpha \pm \beta i \in \mathbb{C}$ for $\beta \in \mathbb{R}$ and $\alpha < 0$ (green), $\alpha = 0$ (orange), and $\alpha > 0$ (blue), respectively).  (b) Several branches of the solution to the Lambert $W$ function for $k = -2, -1, 0, 1, 2$ (given in equation \eqref{lambertsols}) for $D_s = 1$ and $\tau_s \in [0.1 \ 3]$.}
	\label{fig_chareq}
\end{figure}

 Since equation \eqref{homodT} solely depends on the delayed state $T_s(t-D_{s})$, we assume its solution to be an exponential equation of the form 
\begin{equation}
T_s(t;\lambda) = ce^{\lambda t},
\label{assumedsoln}
\end{equation}
where $c$ is a scaling factor and $\lambda$ is the eigenvalue \cite{Asl2003,Bellman1963,Ruan2003}. We make this assumption since the DDE can reduce to an ordinary differential equation (ODE) in steady-state. Hence, we can find an explicit solution to equation \eqref{homodH} as 
\begin{equation}
H(t;\lambda) = \dfrac{H_I H_s}{\tau_H\lambda + 1} T_s(t;\lambda) + \bigg(H_0 - \dfrac{H_I H_s}{\tau_H\lambda + 1} T_{s,0} \bigg) e^{-\frac{1}{\tau_H} t}. 
\end{equation}
Therefore, $H$ is linearly related to $T_s$. As $t \to \infty$, the exponential term vanishes, that is, for large $t$, $H$ is proportional to $T_s$ and the stability of $H$ depends explicitly on the stability of $T_s$. Thus, by analyzing $T_s$, we inherently know the behavior of $H$. By substituting equation \eqref{assumedsoln} into equation \eqref{homodT}, we obtain 
\begin{equation}
\lambda c e^{\lambda t} = -\frac{1}{\tau_s} c e^{\lambda(t-D_s)} \quad \Rightarrow \quad c e^{\lambda t} ( \tau_s \lambda + e^{\lambda D_s} ) = 0.
\label{subinT}
\end{equation}
Trivially, if $c=0$, then $T_s(t) = 0$ is a solution to equation \eqref{subinT}. Considering the portion of equation \eqref{subinT} in parentheses, we obtain the transcendental, nonlinear CE
\begin{equation}
\phi (\lambda) = \tau_s \lambda + e^{-\lambda D_s} = 0. 
\label{chareq}
\end{equation}
Note that $\tau_s, D_s >0$.  Several curves of the CE are plotted in Figure \ref{fig_chareq} for $D_s = 1$ and $\tau_s$ sampled from the interval $[1.25/\pi, 1.25e]$. $\phi$ can have 2 real roots (gray curve), 1 real root (red curve), or infinitely many complex roots. Examples of solutions with complex roots are also plotted in Figure \ref{fig_chareq}, {\it i.e.}, for $\lambda \in \mathbb{C}$, $\lambda = \alpha \pm \beta i$ for $\beta >0$. When $\alpha < 0$ (green curve), the solutions are stable. For $\alpha = 0$ (orange curve), a limit cycle emerges about the origin. When $\alpha >0$ (blue curve), solutions are unstable. 

	\par The CE (equation \eqref{chareq}) has been analyzed in terms of the Lambert $W$ function for $W: \mathbb{C} \rightarrow \mathbb{C}$ \cite{Asl2003,Kuang1993,Yi2006} that satisfies 
	\begin{equation}
	W(z) e^{W(z)} = z
	\label{Lambert}
	\end{equation}
	for $z \in \mathbb{C}$. To obtain a similar formulation as equation \eqref{Lambert}, we rearrange the characteristic equation as 
	\begin{equation} 
	\tau_s \lambda + e^{-\lambda D_s} = 0 \quad \Rightarrow \quad \lambda e^{\lambda D_s} = \dfrac{-1}{\tau_s},
	\end{equation}
	and multiply both sides by $D_s$, which gives 
	\begin{equation}
	\lambda D_s e^{\lambda D_s} = \dfrac{-D_s}{\tau_s}. 
	\label{newchareq}
	\end{equation}
	This equation is a Lambert $W$ function (equation \eqref{Lambert}) for $z = -D_s/\tau_s$, that is,  
	\begin{equation} 
	W\Bigg(-\dfrac{D_s}{\tau_s} \Bigg) e^{W\big(-\frac{D_s}{\tau_s}\big)} = - \dfrac{D_s}{\tau_s}. 
	\label{charLambert}
	\end{equation} 
	Combining equations \eqref{newchareq} and \eqref{charLambert} gives
	\begin{equation} 
	\lambda D_s = W\Bigg(-\dfrac{D_s}{\tau_s}\Bigg) \quad \Rightarrow \quad \lambda = -\dfrac{1}{D_s} W\Bigg(-\dfrac{D_s}{\tau_s}\Bigg). 
	\label{lambertsols}
	\end{equation}
	Note that the Lambert $W$ function contains an infinite number of branches $W_k$ for integer $k = - \infty, \dots, -1, 0, 1, \dots \infty$ with corresponding solutions $\lambda_k$. Figure \ref{lambertsolsfig} shows several solutions to equation \eqref{lambertsols} for $D_s = 1$, $\tau_s$ varied from $[0.1, \ 3]$, and $k = -2, -1, 0, 1, 2$. Since the focus of this study is the analysis of model behavior in relation to patient data, we are concerned with solutions that have the greatest real part, corresponding to the principal branch, $W_0(z)$ \cite{Shinozaki2006}. Therefore, we only analyze the solution for the principal branch, that is, (i) for $\text{Re}(W_0(z)) < 0$, the solutions corresponding to $W_k(z)$ for $k \neq 0$ vanish faster than the solution corresponding to $W_0(z)$; (ii) for $\text{Re}(W_0(z)) = 0$, the solutions corresponding to $W_k(z)$ are either purely oscillatory (as in $W_0(z)$) or vanish; or (iii) for $\text{Re}(W_0(z)) > 0$, the solution corresponding to $W_0(z)$ dominates and diverges faster than all other solutions. \\

\par \noindent \textbf{Case 1: Real Roots.} We consider the case where the CE ( equation \eqref{chareq}) has real roots, {\it i.e.}, $\lambda = \lambda_1, \lambda_2 \in \mathbb{R}$ and $\lambda_1 \neq \lambda_2$. Taking the derivative of $\phi$ and setting it to zero, we obtain, 
\begin{equation} 
\phi'(\lambda) = \tau_s - D_s e^{\lambda D_s} = 0 \quad \Rightarrow \quad \lambda^* = -\dfrac{1}{D_s} \ln\Bigg(\dfrac{\tau_s}{D_s}\Bigg)
\end{equation}
for minimizer $\lambda^*$. Therefore, $\phi$ has real roots if and only if the local minimum of $\phi$ is less than or equal to zero, that is, $\phi(\lambda^*) \leq 0 \Leftrightarrow \phi$ has real roots. Substituting $\lambda^*$ into equation \eqref{chareq} and setting $\lambda^* \leq 0$ yields 
\begin{equation} 
\phi(\lambda^*) = \dfrac{\tau_s}{D_s} \Bigg( -\ln \Big(\dfrac{\tau_s}{D_s}\Big) + 1 \Bigg) \leq 0.
\end{equation} 
Since $\tau_s, D_s > 0$, we have $-\ln(\tau_s/D_s) + 1 \leq 0$ and 	
\begin{equation} 
eD_s \leq \tau_s. 
\label{sinkline}
\end{equation} 
This is analogous to the principal branch of the Lambert $W$ function. Note that $W_0$ is real- and single-valued for arguments greater than or equal to $-1/e$, {\it i.e.},
	\begin{equation}
	-\dfrac{D_s}{\tau_s} \geq -\dfrac{1}{e} \quad \Rightarrow \quad eD_s \leq \tau_s. 
	\end{equation}
When equation \eqref{sinkline} is an equality, $\phi(\lambda)$ has one real root and solutions to the origin are critically damped, analogous to the harmonic oscillator. Theoretically, we observe a transcritical bifurcation at this line in the $D_s \ \tau_s$-plane (Figure \ref{homoContour}), where the critical point changes stability. When equation \eqref{sinkline} is strictly greater than, there are two real solutions to $\phi$ that produce overdamped behavior. Hence, the solutions that obey the constraint given in equation \eqref{sinkline} are stable. \\

\par \noindent \textbf{Case 2: Complex Roots.} We consider the case where the CE has complex roots, that is, $\lambda = \alpha \pm \beta i$ for $\alpha, \beta \in \mathbb{R}$ and $\beta > 0$. Without loss of generality, we consider $\lambda = \alpha + \beta i$. Then,
\begin{align}
0 &= \phi(\lambda) =  \Big(\alpha \tau_s + e^{-D_s\alpha}\cos(D_s\beta) \Big) + i \Big( \beta \tau_s - e^{-D_s\alpha}\sin(D_s\beta) \Big).
\end{align}
Since $\phi(\lambda) = 0$, both the real and imaginary parts of $\phi(\lambda)$ must also equal to $0$. Thus, 
\begin{align}
0 &= \mathrm{Re}(\phi(\lambda)) = \alpha \tau_s + e^{-D_s\alpha}\cos(D_s \beta)
\label{realpart} \ \text{and} \\
0 &= \mathrm{Im} (\phi(\lambda)) = \beta \tau_s - e^{-D_s \alpha} \sin(D_s \beta).
\label{imaginarypart}
\end{align}
Dividing equation \eqref{realpart} by equation \eqref{imaginarypart}, yields
\begin{equation}
\alpha = -\beta \cot(D_s\beta).
\label{complexrelation}
\end{equation}
Therefore, it is guaranteed that $\alpha<0$ as long as $0 + k\pi < D_s\beta \leq \pi/2 + k\pi$ for $k \in \mathbb{Z}$. Since we are concerned with the principal branch $W_0$, we set $k = 0$. Hence, $\alpha < 0$ for $0 < D_s \beta \leq \pi/2$, and the solutions are asymptotically stable and spiral towards the origin. 

\par
\begin{figure}[t!]
	\centering 
	\includegraphics[width=.8\textwidth]{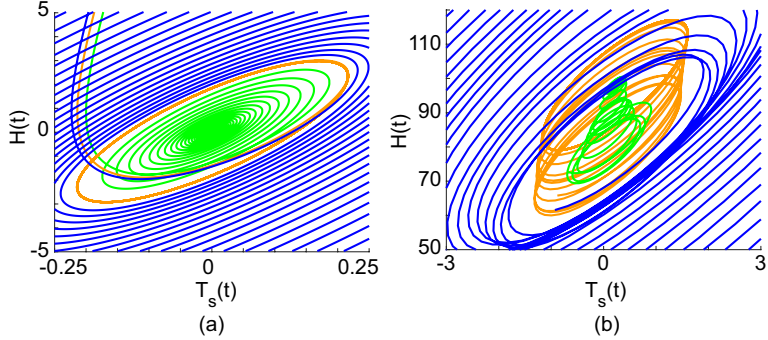}
	\caption{Hopf bifurcation observed in (a) the homogeneous system in equation \eqref{homogeneous} and (b) the nonhomogeneous system in equation \eqref{matsys}. The limit cycle (orange) is plotted along with solutions spiraling out from the critical point (blue) and into the critical point (green). Parameter values for each solution: unstable - $D_s = 1, \tau_s = 1.95/\pi$, limit cycle - $D_s = 1, \tau_s = 2/\pi$, stable focus - $D_s = 1, \tau_s = 2.1/\pi$.}
	\label{hopfBifur}
\end{figure}

 When $\alpha = 0$,  we have from equation \eqref{realpart}
\begin{equation}
\cos (D_s \beta ) = 0 \quad \Rightarrow \quad 
D_s \beta= \frac{\pi}{2} + k\pi \quad \textnormal{ for } k \in \mathbb{Z}.
\label{Dbeta}
\end{equation}
Once again, we set $k = 0$, since we are concerned with the principal branch $W_0$. Then,  $D_{s} \beta = \pi/2$ and characterizes the division between stable and unstable behavior \cite{Asl2003}. Substituting equation \eqref{Dbeta} into equation \eqref{imaginarypart} yields, 
\begin{equation} 
\frac{\pi}{2 D_s} \tau_s = \sin \Big( \dfrac{\pi}{2} \Big) = 1 \quad \Rightarrow \quad D_s = \frac{\pi}{2} \tau_s.
\label{limitcycleline}
\end{equation}
This line is where $\lambda$ crosses the imaginary axis, resulting in a limit cycle. Therefore, theoretically we observe a Hopf bifurcation about the origin in the states (Figure \ref{hopfBifur}). As shown in previous studies \cite{Asl2003,Kuang1993}, instability in the principal branch occurs when 
\begin{equation} 
	D_s > \dfrac{\pi}{2} \tau_s. 
\end{equation}

\par
\begin{figure}[!th]
	\centering 
	\includegraphics[width=.5\textwidth]{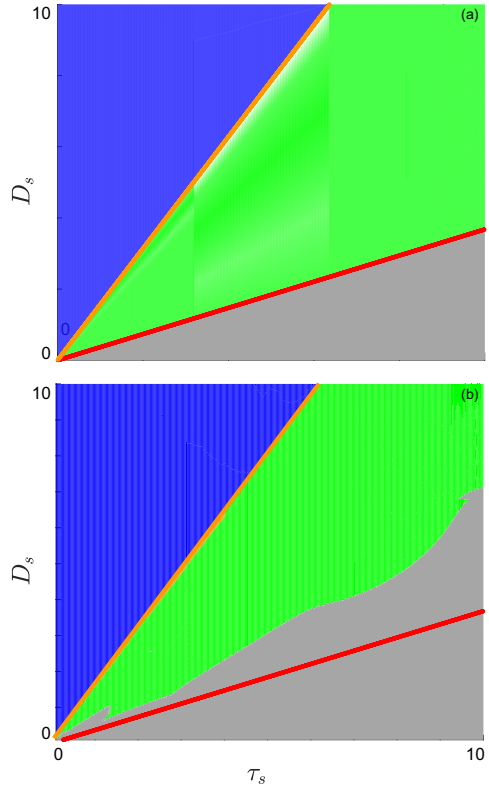}
	\caption{Bifurcation diagrams of the behavior of $T_s$ for $D_s \in [0.1,10]$ and $\tau_s \in [0.1,10]$ evaluating the (a) homogeneous system \eqref{homogeneous} and (b) nonhomogeneous system \eqref{matsys}. Solutions types are denoted as overdamped (gray), critically damped (red line), stable focus (green), limit cycle (orange line), and unstable (blue). The red line indicated in panel (b) denotes the analytically derived line $\tau_s = eD_s$ for comparison to show the increased sink region (gray).}
	\label{contours}
\end{figure}

 In conclusion, we have shown that for the homogeneous system (equation \eqref{homogeneous}) there exist two lines across which the behavior of the solutions changes: one at $eD_s = \tau_s$ where the solutions to the CE change from real to imaginary and one at $D_s = \tau_s \pi/2$ where the imaginary roots cross the imaginary axis. These lines are shown in Figure \ref{contours}a. The former results in a transcritical bifurcation, changing the origin (which is the critical point) from a sink to a stable focus, and the latter results in a Hopf bifurcation, producing a limit cycle about the origin, as shown in Figure \ref{hopfBifur}a. We classify each solution as one of the following types with the corresponding color from the contour in Figure \ref{contours}a: 
\begin{itemize}
	\item {\it Sink - overdamped} (gray): When $e D_s < \tau_s$, $\phi$ has two real solutions $\lambda_1, \lambda_2 < 0$.
	\item {\it Sink - critically damped} (red): When $e D_s = \tau_s$, $\phi$ has one real solution $\lambda <0$, determining a transcritical bifurcation. 
	\item {\it Stable focus} (green): When $\tau_s/e < D_s < \tau_s \pi/2$,	$\phi$ has complex solutions $\lambda = \alpha \pm \beta i$ and $\alpha < 0$.
	\item {\it Limit cycle} (orange): When $D_s = \tau_s \pi/2$, $\phi$ has complex solutions $\lambda = \pm \beta i$ ({\it i.e.}, $\alpha = 0$) and a limit cycle about the critical point emerges from a Hopf bifurcation.
	\item {\it Unstable} (blue): When $D_s > \tau_s \pi/2 $, $\phi$ has complex solutions $\lambda = \alpha \pm \beta i$ for $\alpha > 0$ and solutions diverge. 
\end{itemize}

\subsection{Nonhomogeneous system}

\par \noindent The inclusion of the forcing function, $\mathbf{f}$, complicates the analysis, and as discussed in Section \ref{Intro}, tools such as \texttt{DDE-Biftool} \cite{Engelborghs2002} and \texttt{knut} \cite{Szalai2006} are not suitable for the nonautonomous, stiff system given in equation \eqref{matsys}. As shown in Figure \ref{forcingfcns}, $\mathbf{f}$ relies on a polynomial fitted to blood pressure data with baseline extended before and after the dynamic behavior. This forcing function creates a set point at a prescribed equilibrium calculated from the baseline SBP ($\bar{P}$) and heart rate ($\bar{H}$). A disturbance caused by this function, such as the VM, can result in undesirable model behavior and instability. Moreover, the perturbation of this control system may cause {\em persistent instability}, that is, instability as a result of a perturbation caused by a forcing function in which oscillatory behavior arises that either remains unstable, oscillates with constant amplitude, or takes a long time to dampen in relation to the stimulus.

\par In the previous section, we determined regions of the parameter space where the five different behaviors arise for the homogeneous solution. We do not expect these regions to be the same for the nonhomogeneous solutions, especially since forcing functions can stabilize and destabilize systems \cite{Sipahi2011}. However, we do expect analogous regions corresponding to the behaviors given above.  

\par This numerical experiment explores the impact of sudden, transient effects of the forcing function on the stability of the system in equation \eqref{matsys} given a specified parameter range. For this analysis we consider only the effects of the interactions between $D_s$ and $\tau_s$. We chose these parameters to investigate based on the analysis of the homogeneous system, which created stability subregions in the parameter space (equations \eqref{sinkline} and \eqref{limitcycleline}). We assume changing these parameters will also cause instability in the nonhomogeneous system. The parameter space for $\tau_s$ and $D_s$ is $[0.1,10] \times [0.1, 10]$ with a discretized mesh of step-size $h = 10^{-3}$ with all other parameters remaining constant at their nominal values. The model was evaluated iteratively at every point in the mesh. Since the instability is most prominent after the maneuver and during the recovery, we only consider the stability of the signal after the breath hold of the VM was over. 

\par Numerically, we have developed an algorithm to determine the type of solution behavior summarized in Algorithm \ref{alg1}. Of particular note are the thresholds $\eta_1$, $\eta_2$, and $\mu$. $\eta_1 = 0.5$ and $\eta_2 = -10^{-2}$ are the maximum and minimum thresholds for the slope of the regression line to determine a limit cycle. $\mu = 0.8$ is the threshold for the $r^2$ value of the line of regression of the amplitudes of consecutive oscillations determining the goodness of fit. These thresholds ensure a limit cycle is obtained. 

\par Solutions for the homogeneous and nonhomogeneous equations were calculated using the stiff, delay differential equation solver RADAR5 \cite{Guglielmi2001}. This is a variable-step solver that employs collocation methods to calculate the history of the delayed states. All initial conditions and constant history value (Table \ref{init}) were assigned such that the system begins in steady-state.

\begin{figure}[!t]
	\centering 
	\includegraphics[width=\textwidth]{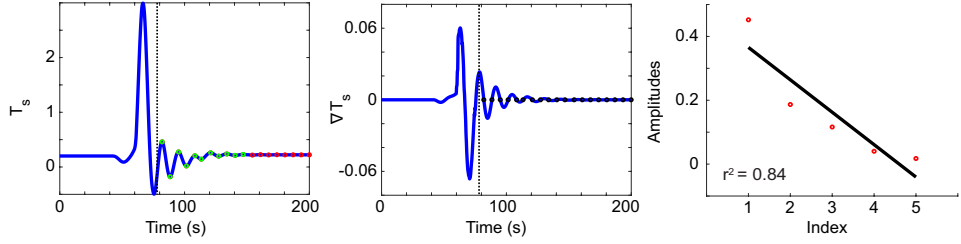}
	\caption{Plots denoting different stages in Algorithm \ref{alg1}. The end of phase III of the Valsalva maneuver is indicated with the vertical dashed line. (a) Representative solution of $T_s$ exhibiting a stable focus with local extrema occurring after the Valsalva maneuver (green circles). The red circle indicates an excluded point if the difference between it and the preceding extremum is $\leq 10^{-8}$. (b) $\nabla T_s$ with zeros indicating local extrema (black circles).  (c) Consecutive amplitudes (red dots) and a line of regression (black line) with the $r^2$ value indicated.}
	\label{algfig} 
\end{figure}

\begin{algorithm}[p]
	\caption{Determine the type of solution behavior for $T_s$.}
	\footnotesize
	\begin{enumerate} 
		\item Calculate $\nabla T_s$ (Figure \ref{algfig}b). Only consider $\nabla T_s$ after the Valsalva maneuver. 
		\item Determine where $\nabla T_s$ crosses the $x$-axis (Figure \ref{algfig}b).
		\item Filter out local extrema if the distance between consecutive points is $<$ 0.1 s. 
		\item \begin{algorithmic}[1]
			\State Let $\mathbf{M}$ and $\mathbf{m}$ be vectors of the local maxima and minima and $N = \text{size}(\mathbf{M})$. 
			\For{$i = 1$ to N}
			\If{ $|M_i - m_i| < 10^{-8}$} 
			\State  Remove $M_i$ and $m_i$ from $\mathbf{M}$ and $\mathbf{m}$, respectively. 
			\EndIf
			\EndFor
			\State $\tilde{\mathbf{M}}$ and $\tilde{\mathbf{m}}$ are the resulting filtered vectors. 
		\end{algorithmic}
		\item  Determine the vector of amplitudes $\mathbf{a} = \tilde{\mathbf{M}} - \tilde{\mathbf{m}}$. 
		\item Assign solution behavior. 
		\begin{algorithmic}[1] 
			\If{$\mathbf{a}$ is empty} 
			\State $T_s$ is a sink. 
			\ElsIf{$\mathbf{a}$ has 1 entry}
			\State $T_s$ spirals in. 
			\Else 
			\State Fit a regression line through the entries of $\mathbf{a}$, $y = b_0 + b_1 x$ (Figure \ref{algfig}c).
			\State Calculate the $r^2$ value of the regression line. 
			\If{ $\eta_2 \leq b_1 \leq \eta_1$} 
			\If{ $r^2 > \mu$} 
			\State $T_s$ is a limit cycle. 
			\Else
			\State $T_s$ spirals in. 
			\EndIf 
			\ElsIf{ $b_1 > \eta_1$} 
			\State $T_s$ spirals out. 
			\Else 
			\State $T_s$ spirals in. 
			\EndIf 
			\EndIf 
		\end{algorithmic} 
	\end{enumerate} 
	\label{alg1} 
\end{algorithm} 

\begin{figure}[!h]
	\centering
	\includegraphics[width=.8\textwidth]{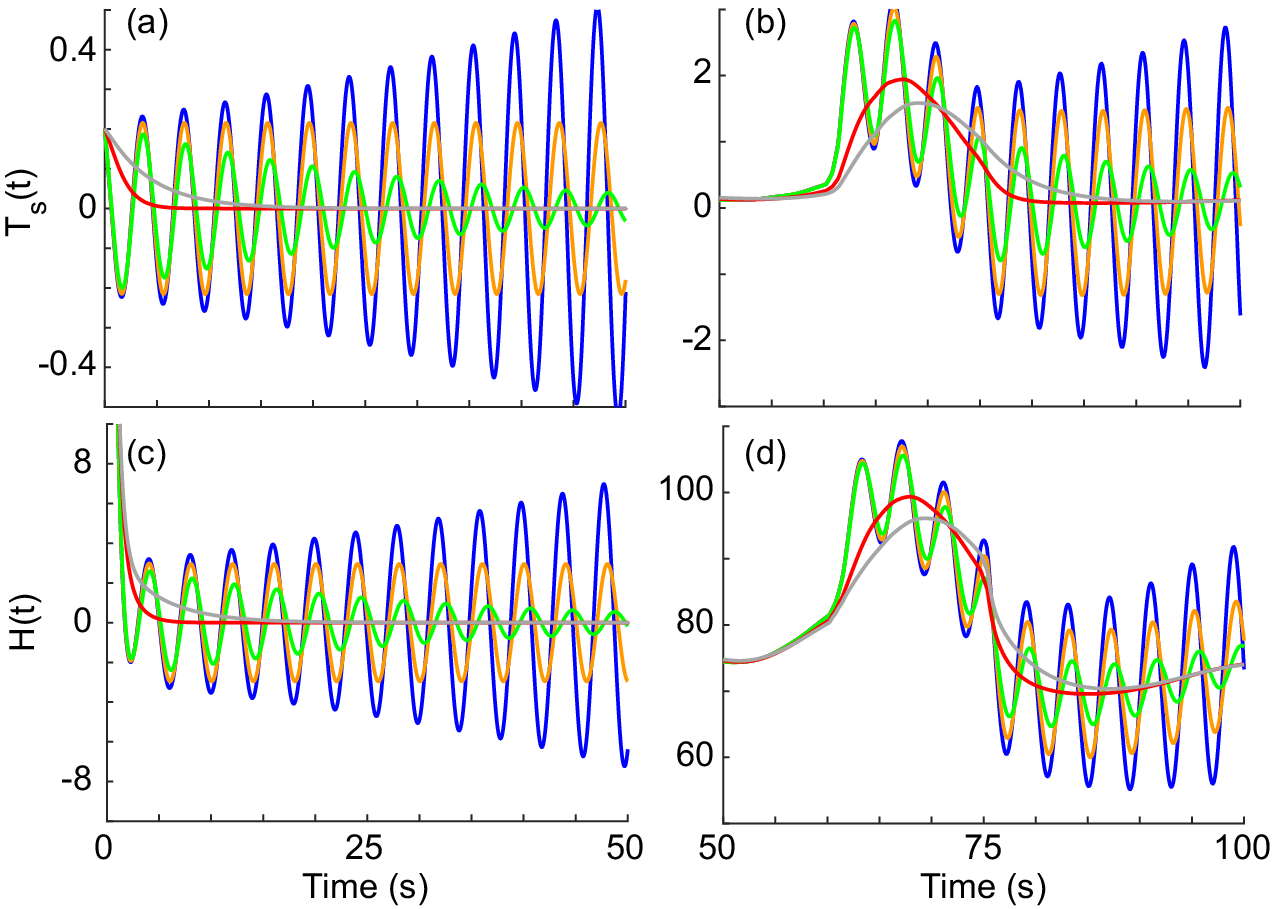}
	\caption{Representative solutions from each of the stability regions given in Figure \ref{contours} with corresponding colors and $D_s = 1$: overdamped sink (gray, $\tau_s = 2e$), critically damped sink (red, $\tau_s = e$), stable focus (green, $\tau_s = 2.1/\pi$), limit cycle (orange, $\tau_s = 2/\pi$), and unstable (blue, $\tau_s = 1.95/\pi$). (a and c) $T_s$ and $H$ for the homogeneous equation \eqref{homogeneous}. (b and d) $T_s$ and $H$ for the nonhomogeneous system \eqref{matsys}.}
	\label{modelsims}
\end{figure} 

\section{Results and discussion} 

\par \noindent Stability analysis is important for understanding the kinds of outputs a model can produce. Moreover, bifurcation analysis explores where changes in parameter values impact the system. In this study, we have analyzed a system of nonautonomous, stiff, delay differential equations (DDEs) that can be written as the sum of a homogeneous system (equation \eqref{homogeneous}) and its associated nonhomogeneous component including the forcing vector $\mathbf{f}$ given in equation \eqref{forcingvector}. Figure \ref{contours} displays the results of the stability analysis with a two-parameter bifurcation diagram plotted for both the homogeneous system (Figure \ref{contours}a) and the nonhomogeneous system (Figure \ref{contours}b in equation \eqref{reducedsystem}) denoted with the following colors: 
\begin{itemize} 
	\item {\it Sink - overdamped} (gray): The inclusion of the forcing function $\mathbf{f}$ increases the range of the stable region within the parameter space, that is, the gray region extends beyond the red line in Figure \ref{contours}b. 
	\item {\it Sink - critically damped} (red): There is a shift in the red line denoting the transcritical bifurcation from the homogeneous contour to the nonhomogeneous contour, as the sink region expands. This relation is no longer a line but a curve between the sink and stable focus regions. 
	\item {\it Stable focus} (green): The stable focus region shrinks with the inclusion of the forcing function, resulting in an oscillatory contour that is steeper than the red line predicted in the homogeneous system. 
	\item {\it Limit cycle} (orange): The limit cycle occurs in the same location in both the homogeneous and the nonhomogeneous bifurcation contours. This is most likely due to the fact that the large amplitude oscillations begin to dominate the signal. 
	\item {\it Unstable} (blue):  The unstable region is the same in both the homogeneous and nonhomogeneous contours. This is to be expected as the solutions diverge. 
\end{itemize}  
Figure \ref{modelsims} displays representative curves from each region mentioned above for both $T_s$ (Figure \ref{modelsims}a and \ref{modelsims}b) and $H$ (Figure \ref{modelsims}c and \ref{modelsims}d). For the nonhomogeneous system (Figures \ref{modelsims}b and \ref{modelsims}d), the solutions begin in steady-state and the forcing function $\mathbf{f}$ induces the VM, causing some of the responses to have oscillatory behavior. Holding $D_s$ constant at its nominal value, we varied $\tau_s$ showing that as $\tau_s$ decreases, we see a shift in the behavior of the model output from sink (gray) to stable focus (green) to limit cycle (orange) to unstable (blue). We observe that though $D_s$ and $\tau_s$ are within their individual physiological ranges, their interactions cause persistent instability for decreasing values of $\tau_s$. Therefore, to ensure that the model produces physiologically relevant results, restricting the parameter space to remain in the sink and stable focus regions is necessary. 

\par Bifurcation packages, {\it e.g.}, \texttt{DDE-Biftool} \cite{Engelborghs2002}, cannot be used to conduct the analysis of the nonautonomous, stiff system of delay differential equations discussed here. Moreover, the forcing function is not periodic, so packages, such as \texttt{knut} \cite{Roose2007,Szalai2006}, cannot be utilized to analyze equation \eqref{reducedsystem}. Therefore, we developed our own algorithm to qualitatively assess the behavior of the solutions propagated after the implementation of the VM. This algorithm classifies the behavior of the solutions starting after the transient VM stimulus by quantifying consecutive amplitudes. With this algorithm, we were able to effectively categorize the solutions and determine the boundary between the sink and stable focus regions.  

\begin{table}[!b]
	\centering 
	\footnotesize  
	\begin{threeparttable} 
		\caption{Parameter values for control subjects and a postural orthostatic tachycardia syndrome (POTS) patient with M behavior.}
		\begin{tabular}{clcccc}
			\hline
			{\bf Subject} & \multicolumn{1}{c}{{\bf Classification}}  & $D_s$ & $\tau_s$ & {\bf Stability region} & {\bf Color}\\ 
			\hline 
			\hline 
			1& Control &  9.2  & 7.5  & Sink & Gray \\ 
			2& Control & 4.7  &  5.4  & Stable focus & Green \\
			3 & POTS with M behavior &  5.6 &  5.2 & Stable focus & Green \\
			\hline 
		\end{tabular} 
		\label{parsresults}
	\end{threeparttable} 
\end{table} 

\begin{figure}[!t] 
	\centering  
	\includegraphics[width=\textwidth]{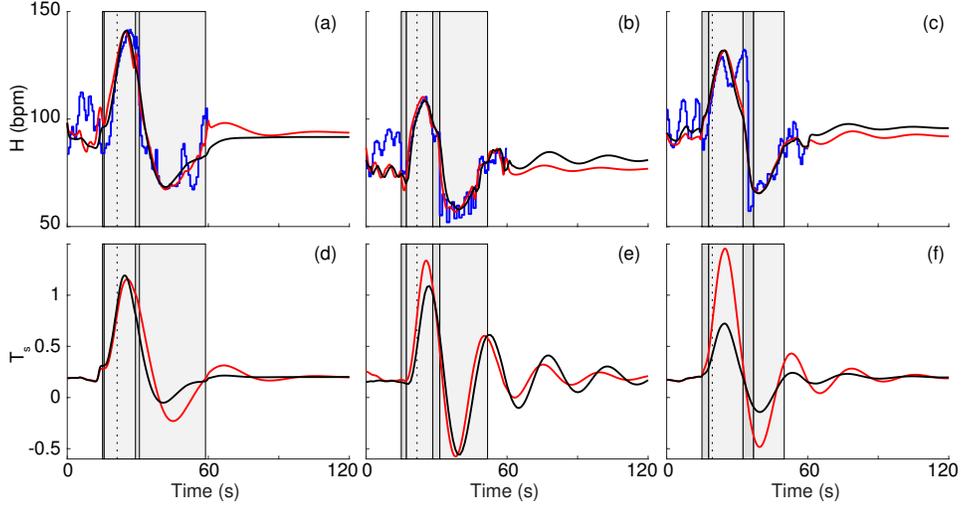}
	\caption{Heart rate ($H$, bpm) model fits for the reduced two-state model (red, equation \eqref{reducedsystem}) and the full five-state model (black, equation \eqref{fullsystem}) to heart rate data (blue) and the resulting $T_s$ trace for (a and d) Subject 1 - control and sink, (b and e) Subject 2 - control and stable focus, and (c and f) Subject 3 - postural orthostatic tachycardia syndrome (POTS) and stable focus. Solutions are calculated using the two-dimensional nonhomogeneous system .} 
	\label{datacomp} 
\end{figure} 

\par Figure \ref{datacomp} shows the heart rate model fits for the reduced two-state model (red, equation \eqref{reducedsystem}) and full five-state model (black, equation \eqref{fullsystem}) and the resulting $T_s$ trace for each subject: Subject 1 is a control subject exhibiting sink behavior (Figures \ref{datacomp}a and d); Subject 2 is a control subject exhibiting stable focus behavior (Figures \ref{datacomp}b and e); and Subject 3 has POTS also exhibiting a stable focus behavior (Figures \ref{datacomp}c and f). The model was solved for 120 s to extend the signal and allow oscillations sufficient time to dampen. Parameter values and where they fall in the bifurcation diagram (Figure \ref{contours})  are given in Table \ref{parsresults}. As shown in Figure \ref{datacomp}, control subjects can have both sink and stable focus behavior. This could be due to the fact that some subjects have naturally higher baseline sympathetic activity or due to undiagnosed autonomic dysfunction.

\par The M behavior is hypothesized to be the result of overactive sympathetic and parasympathetic activity \cite{Palamarchuk2016}. Figure \ref{datacomp}f demonstrates that the baroreflex control mechanism is very sensitive for this subject, causing oscillatory dynamics that are generally not seen in most control subjects \cite{Randall2019}. With the stability analysis, we are able to characterize the dynamics of the M behavior as well as support that the M behavior may be due to instability in the negative feedback control of the baroreflex stimulating a sympathetic response. 

\par When comparing the full and reduced model outputs in Figure \ref{datacomp}, we observe that the full model induces more smoothing of the heart rate than the reduced model. The amplitudes of the reduced model are larger than those of the full model. This is to be expected, since the full model has more differential equations with varying time-scales, which can smooth the solutions. Despite this, the reduced model displays a good approximation of the full model. The parameter values used in the both model formulations are given in Table \ref{parsresults}. Note that the full and reduced models agree in the type of signal behavior. Hence, the reduced model with its inherent simplicity can successfully detect persistent instability in these subjects. 

\par We do not see limit cycles or unstable modes in practice. This is most likely due to the fact that when one system becomes inordinately overactive, there are many other redundancies in place to reset the body, such as inducing syncope \cite{Boron2017}. However, we can classify each of these stability regions not only on the basis of their mathematical properties but of their clinical relevance. These categories are: 
\begin{itemize}
	\item {\it Sink}: Healthy/control behavior within the ``normal" range. 
	\item {\it Stable focus}: Potential dysfunction caused by overactive sympathetic behavior.
	\item {\it Limit cycle}: Unphysiological steady pulsation of sympathetic activity. 
	\item {\it Unstable}: Sympathetic positive feedback that may be unphysiological. (If it is physiological, it may be corrected via other regulatory mechanisms or, in a worst case scenario, cause death). 
\end{itemize} 

\par In this study, we physiologically reduced the full five-state model to a system of two DDEs that can be solved analytically. This simplified both the model and the analysis and proved to be a reasonable reduction of the system (Figure \ref{fullvstwostate}). However, numerically, we could have analyzed the full model as opposed to the reduced two-state model. We considered the two-state model for both the homogeneous and nonhomogeneous analyses to facilitate interpretability and comparison. 

\section{Conclusions} 

\par \noindent In this study, we analyzed the effect of the delay differential equation system modeling the autonomic response to the Valsalva maneuver and categorized the various types of behavior that can result from the interaction of the delay parameter $D_s$ and the time-scale $\tau_s$. Moreover, we classified stability regions both mathematically and physiologically in a two-parameter bifurcation diagram. Motivated by oscillatory behavior that arises in the data, we have determined a numerical relationship between $D_s$ and $\tau_s$ and observed transcritical and Hopf bifurcations. The model also supports that the M behavior may arise due to oscillatory behavior from overactive sympathetic stimulation. 

\appendix

\section*{Appendix}

\setcounter{equation}{0}
\renewcommand{\theequation}{A\arabic{equation}}

\par \noindent Table \ref{parameters} lists the nominal parameter values used in the full five-state model. Parameters $A$ (dimensionless), $q_w$ (mmHg$^{-1}$), and $s_w$ (mmHg) are associated with equation \eqref{epswj} and informed by the sigmoid-like relationship given in Mahdi {\it et al.} \cite{Mahdi2013}. However, in our formulation, $s_w$ is a half-saturation value defined as the mean systolic blood pressure ($\bar{P}$, mmHg), which changes for each subject as listed in Table \ref{datatable}. $B$ (s$^{-1}$, equation \eqref{n}) varies from 0 to 1, denoting the convex combination of the two afferent signals. We set $B = 0.5$ to model average signaling. Nominal values for the time-scales of the afferent and efferent neural signals were set to literature values: $\tau_b = 0.9$, $\tau_{p} = 1.8$, and $\tau_s = 10$ s \cite{Lu2001}. $H_I = 100$ bpm \cite{Olufsen2005} is the intrinsic heart rate, {\it i.e.}, when the heart is completely denervated. The sigmoid steepness parameters for the parasympathetic and sympathetic neural signals ($q_p$ (s) and $q_s$ (s)) were set to 10 s to account for the order of magnitude of $n$. The sympathetic delay ($D_s$) was assigned a value of 3 s as used in previous modeling studies \cite{Lu2001,Wesseling1993}

\par Some parameter values were calculated a priori given baseline values ($\bar{P}$ (mmHg) and $\bar{H}$ (bpm)) from the data. From here on, $\bar{\cdot}$ denotes a steady-state value. Using this notation, the following parameters are calculated as: 
\begin{align}
s_p &= \bar{n} + \ln(K_p / T_{p,0} - 1)/q_p, \label{s_p} \\
s_s &= \bar{n} - \ln(K_s / T_{s,0} - 1)/q_s, \label{s_s}\\ 
H_p &= (1 - \bar{H}/H_I + H_s T_{s,0})/T_{p,0}, \quad \text{and} \label{H_p}\\ 
H_s &= (H_M/H_I - 1)/K_s, \label{H_s} 
\end{align}
where $T_{p,0}$ and $T_{s,0}$ are the initial condition and constant history value given in Table \ref{init} and $H_M$ (bpm) is the maximal heart rate as a function of age \cite{Tanaka2001}. Further information on parameter value assignment can be found in Randall {\it et al.} \cite{Randall2019}.

\section*{Acknowledgments}

\par \noindent We thank Drs. Jesper Mehlsen, Section of Surgical Pathophysiology, Juliane Marie Centre, Rigshospitalet, University of Copenhagen, Denmark, and Louise S. Brinth, Department of Clinical Physiology and Nuclear Medicine, Bispebjerg and Frederiksberg Hospital, Denmark for providing the data analyzed in this study. 

\section*{Funding}

\par \noindent This study was supported by the Division of Mathematical Sciences at the National Science Foundation grant \#1246991.



\bibliographystyle{elsarticle-num} 


\end{document}